\documentclass[preprint,11pt]{elsarticle}
\usepackage[latin1]{inputenc}
\usepackage{amsthm}
\usepackage{amssymb}
\usepackage[english,frenchb]{babel}
\usepackage[left=1.5cm,right=1.5cm,top=4cm,bottom=4cm]{geometry}
\input xy
\xyoption{all}
\newtheorem{thm}{Théorème}
\newtheorem{hyp}[thm]{Hypothèses}
\newtheorem{lemme}[thm]{Lemme}

\newtheorem{definition}[thm]{Définition}
\newtheorem{remarque}[thm]{Remarque}

\begin{document}
\begin{frontmatter}
\title{Actions modérées de schémas en groupes affines et champs modérés \\ Tame actions of affine group schemes and tame stacks}
\author{Sophie MARQUES}
\address[MARQUES]{Univ. Bordeaux, IMB, UMR 5251, F-33400 Talence, France et Università degli studi di Padova, Dipartimento di Matematica, Via Trieste 63, 35121, Padova, Italia dans le cadre de ALGANT.DOC. }
\ead{Sophie.Marques@math.u-bordeaux1.fr}
\begin{abstract}
We compare two notions of tameness: one introduced in \cite{Boas} for actions of affine group schemes and one introduced in \cite{AOV} for stacks. From this comparison, we deduce results on the structure of inertia groups in these tame situations.\\
Nous comparons ici deux notions de modération : celle d'action modérée introduite dans \cite{Boas} et celle de champ modéré introduite dans \cite{AOV}, pour ensuite en déduire des résultats de structure sur les groupes d'inertie. 
\end{abstract}
\end{frontmatter}
\setcounter{page}{1} 
\setcounter{tocdepth}{4}
\bibliographystyle{alpha-fr}
\pagestyle{plain}
{\small
 \section*{Remerciements}
Cette note doit énormément à l'aide précieuse et aux conseils de Boas Erez, Marco Garuti et Jean Gillibert. J'ai aussi bénéficié de remarques pertinentes de Brian Conrad, Cédric Pépin et Matthieu Romagny. C'est pourquoi je les remercie vivement ainsi que tous ceux qui ont contribué à l'élaboration de cette note.}
\section{Introduction} 
Chinburg, Erez, Pappas et Taylor ont défini dans leur article \cite{Boas} la notion d'action modérée. De leur côté, Abramovich, Olsson et Vistoli ont introduit dans leur article \cite{AOV} la notion de champ modéré. Il est alors naturel de comparer ces deux notions de modération. En se plaçant sous les hypothèses de \cite[Definition 3.1]{AOV} pour le champ quotient associé à une action donnée, on montre que si l'action est modérée au sens de \cite{Boas}, alors le champ quotient est modéré au sens \cite{AOV}. Sous des hypothèses plus restrictives mais souvent suffisantes dans la pratique, notamment de finitude sur le schéma en groupes qui agit, on montre qu'il y a même équivalence entre ces deux notions. Ces résultats nous permettent ensuite d'établir un théorème de structure pour les groupes d'inertie sous des hypothèses de modération. 
\section{Notations}
Dans la suite, les schémas considérés seront affines sur une base affine. Plus précisément, la base sera $S:= Spec (R)$ où $R$ est un anneau commutatif unitaire, $G=Spec (A)$ sera un schéma en groupes affine plat sur $S$ et $X:= Spec (B)$ sera un schéma affine sur $S$. La donnée d'une action de $G$ sur $X$ définie par un morphisme $\mu_X : X \times_S G \rightarrow X$ est équivalente à la donnée d'un $B$-comodule défini par le morphisme structural, que l'on notera $\rho_B : B \rightarrow B \otimes_R A$. On notera cette action par $(X, \mu_X)$. On notera $C:=B^A :=\{ b \in B |\rho_B ( b) = b \otimes 1\}$ l'anneau des invariants pour l'action et $Y := Spec (C)$.  Si $S'$ est un $S$-schéma, on notera $X_{(S')}= X \times_S S'$. Enfin, l'unité d'un anneau $\Lambda$ sera notée par $1_{\Lambda}$. 
\section{Définitions et résultats}
On rappelle qu'un morphisme $\alpha :(B, \rho_B) \rightarrow (C , \rho_C)$ de $A$-comodules est une application $R$-linéaire telle que $\rho_C \circ \alpha = (\alpha \otimes 1_B) \circ \rho_B$. La $R$-algèbre $A$ peut être vu comme un $A$-comodule via la comultiplication $\Delta : A \rightarrow A \otimes_R A $.
\begin{definition} 
On dira qu'une \textit{action $(X, \mu_X)$ est modérée} s'il existe un morphisme de $A$-comodules $\alpha : A \rightarrow B$, qui est unitaire, \textrm{i.e.} $\alpha (1_A ) = 1_B$. 
\end{definition} 
Dans le cas d'un groupe constant, cette définition peut se traduire en termes de surjectivité du morphisme trace, ce qui étend le critère de Noether caractérisant les extensions d'anneaux d'entiers modérées (cf. \cite[Introduction]{Boas}) et explique ainsi le choix de la terminologie.\\ 
\indent Considérons le champ quotient $[X/G]$ associé à l'action $(X, \mu_X)$. On rappelle qu'un \textit{espace de modules grossier} pour $[X/G]$ est un couple $(M, \rho)$ où $M$ est un espace algébrique et $\rho : [X/G] \rightarrow M$ est un morphisme universel pour les morphismes de $[X/G]$ vers un espace algébrique tel que pour tout corps $\Omega$ algébriquement clos, $|[X/G]( \Omega )| \simeq M (\Omega )$ où $|[X/G](\Omega )|$ est l'ensemble des classes d'isomorphismes du groupoïde. Dans la suite, on notera simplement $\rho$ cet espace de module. Rappelons aussi la définition de groupe d'inertie. Pour $\zeta$ un $T$-point de $X$ où $T$ est une $R$-algèbre, on notera $I_G (\zeta)$ le groupe d'inertie de l'action au point $\zeta : Spec (T) \rightarrow X$, défini comme le produit fibré $$I_G (\zeta):=(X\times_SG )\times_{(X \times_S X , (\mu_X , p_1), (\zeta , \zeta ) \circ \Delta)} Spec (T)$$ où $p_1$ est la première projection.

\begin{hyp} 
\label{hyp} Supposons que le champ quotient $[X/G] $ soit un champ algébrique localement de présentation finie et que tous les groupes d'inertie soient finis. 
\end{hyp}
Sous ces hypothèses, on sait par l'article \cite{KM} qu'il existe un espace de modules grossier $\rho : [X/G]\rightarrow M$ et que le morphisme $\rho$ est propre.
\begin{definition} 
\label{AOVdef}
Sous les hypothèses \ref{hyp}, on dit que le \textit{champ $[X/G]$ est modéré} si le foncteur entre les catégories de faisceaux quasi-cohérents $\rho_* : Qcoh [X/G] \rightarrow Qcoh (M)$ est exact. En particulier, si on considère l'action triviale avec $G$ fini, localement libre sur $S$, on dira que $G$ est \textit{linéairement réductif au sens de \cite{AOV}}, si le champ classifiant $[S/G]$ est modéré.
\end{definition}
\begin{remarque} 
\label{AOVrk}
Par le critère d'Artin  \cite[Théorème 10.1]{Laumon}, on montre que lorsque $G$ est plat, de présentation finie sur $S$, $[X/G]$ est un champ algébrique, on a même que le morphisme canonique $p : X \rightarrow [X/G]$ est un $G$-torseur quasi-compact et représentable. De plus, si l'on suppose $X$ de présentation finie sur $S$, le champ quotient l'est aussi. 
\end{remarque}
Le but de cet article est d'abord de relier ces deux notions de modération pour ensuite en déduire des résultats sur les groupes d'inertie. 
\begin{thm}
\label{thm1}
Supposons que $S$ soit noethérien, que le schéma en groupes affine $G$ soit plat, de type fini sur $S$, que $X$ soit de type fini sur $S$ et que tous les groupes d'inertie soient finis. Si l'action $(X, \mu_X)$ est modérée alors le morphisme $[X/G] \rightarrow Y$ est un espace de modules grossier pour $[X/G]$ et le champ quotient $[X/G]$ est modéré. 
\end{thm} 
Sous des conditions plus restrictives, il y a même l'équivalence entre les deux notions.
\begin{thm} \label{thm2}
Supposons que $Y$ soit noethérien, que le schéma en groupes affine $G$ soit \emph{fini, localement libre} sur $S$ et que le morphisme $X \rightarrow Y$ soit plat. Sous ces conditions, le morphisme $[X/G] \rightarrow Y$ est un espace de modules grossier et l'action est modérée si et seulement si le champ quotient est modéré. 
\end{thm}
\begin{remarque}\begin{enumerate} \item Le sens direct est vrai sans l'hypothèse $X\rightarrow Y$ plat et $Y$ noethérien. 
\item On peut remplacer l'hypothèse $Y$ noethérien par $X$ de type fini si l'on suppose la base $S$ noethérienne. En effet, par \cite[Theorem 3.1 (2)]{conrad}, $Y$ est alors de type fini sur $S$ donc noethérien.
\end{enumerate}
\end{remarque}
En utilisant l'article \cite{AOV}, on démontre le prochain théorème motivé par l'observation suivante. Si l'on considère le cas du schéma en groupes constant associé à un groupe fini $\Gamma$ tel que $B$ et $C= B^\Gamma$ sont des anneaux de Dedekind, on peut montrer que les groupes d'inertie aux points fermés de $X$ sont linéairement réductifs si et seulement si l'extension $B/C$ est modérée.
\begin{thm}\label{thm3}
Supposons que $S$ soit noethérien, le schéma en groupes $G$ soit plat, de type fini sur $S$, que $X$ soit de type fini sur $S$ et que tous les groupes d'inertie soient finis et plats. Le champ quotient est modéré si et seulement si tous les groupes d'inertie sont linéairement réductifs. 
\end{thm}
\begin{remarque}
\begin{enumerate}
\item On peut montrer qu'une action par un schéma en groupes diagonalisable est toujours modérée, et en appliquant \cite[exposé IX , \S 8]{sga}, on montre que ses groupes d'inertie en tout point fermé sont diagonalisables. Cet exemple pointe vers une généralisation possible de ce dernier théorème. 
\item Les théorèmes exposés ici peuvent se généraliser au cas non affine (voir \cite[Definition 7.1]{Boas}).
\end{enumerate}
\end{remarque}
\section{Preuve du théorème \ref{thm1}}
On commence par deux lemmes.
\begin{lemme} \label{lem1}
L'exactitude du foncteur $\rho_* : Qcoh ([X/G]) \rightarrow Qcoh (Spec (B^A))$ est équivalente à celle du foncteur des invariants $(-)^A : B$-$A$-modules $\rightarrow B^A$-modules (pour les notations voir \cite[\S 2]{Boas}).
\end{lemme} 
\begin{proof}Elle est conséquence du fait que $Qcoh ([X/G]) \simeq Qcoh^G (X)$ par \cite[Example 7.17]{Vistoli} et \cite[\S 5]{Hartshorne}.\end{proof}
Le lemme qui suit est plus général que nécessaire mais intéressant en tant que tel (on ne fait pas d'hypothèse de noethérianité sur la base $S$). 
\begin{lemme} 
Supposons le groupe $G$ de présentation finie et plat sur $S$ et $(X, \mu_X)$ une action modérée, alors \\ $ \rho_* : Qcoh ([X/G]) \rightarrow Qcoh (Y)$ est un foncteur exact et $\mathcal{O}_Y\simeq \rho_* \mathcal{O}_{[X/G]}$. \end{lemme} 
\begin{proof}Le résultat est conséquence du Lemme \ref{lem1} au vu de \cite[lemma 2.3]{Boas}.  \end{proof}
Alper appelle un espace de modules ayant la propriété du lemme \textit{un bon espace de modules} (voir \cite{Alper}).\\
Nous pouvons donc utiliser \cite[Theorem 6.6]{Alper}, qui montre l'universalité d'un bon espace de modules pour les morphismes de $[X/G]$ vers un espace algébrique, vu que $[X/G]$ est  noethérien, puisque $X$ l'est. Par suite, $[X/G] \rightarrow Spec (C)$ est universel pour les morphismes de $[X/G]$ vers un espace algébrique. D'ailleurs, sous les conditions de notre théorème \ref{thm1}, on peut appliquer \cite{KM}, qui montre l'existence d'un espace de modules grossier. Par universalité pour les morphismes de $[X/G]$ vers un espace algébrique, il est égal à $Spec(C)$. Ce qui montre que $[X/G]\rightarrow Spec (C)$ est un espace de modules grossier. La suite découle aisément du Lemme \ref{lem1} combiné avec \cite[Lemma 2.3]{Boas}.
\section{Preuve du théorème \ref{thm2}}
Le sens direct se montre comme pour le théorème \ref{thm1} sans l'hypothèse $X\rightarrow Y$ plat et $Y$ noethérien, une fois que l'on a prouvé que $[X/G] \rightarrow Spec(C)$ est un espace de modules grossier. Ceci résulte de \cite[\S 3]{conrad}. La preuve de la réciproque se base sur le lemme algébrique suivant, étant donné le Lemme \ref{lem1}.
\begin{lemme} \label{relativemt} \label{corexact}
Soit $(X, G)$ une action sur $S$. Supposons $G$ fini, localement libre, $C$ noethérien et $B$ plat sur $C$. Alors les assertions suivantes sont équivalentes : 
\begin{enumerate}
\item[(1)] L'action $(X, \mu_X)$ est modérée.
\item[(2)] $B$ est $A$-coplat en tant que $R$-module (cf. \cite[10.8]{wisbauer}). 
\item[(3)] $(-)^A : {}_B \mathcal{M}^A\rightarrow {}_C \mathcal{M}$, $N \mapsto (N)^A$ est exact .
\item[(1')] L'action $(X , G_Y)$ est modérée. 
\item[(2')] $B$ est $A_C$-coplat en tant que $C$-module. 
\item[(3')] $(-)^{A_C} : {}_B \mathcal{M}^{A_C}\rightarrow {}_C \mathcal{M}$, $N \mapsto (N)^{A_C}$ est exact.
\end{enumerate}
\end{lemme}
\begin{proof}
Pour $(1) \Rightarrow (2)$, on applique \cite[Theorem 1.6]{doiI} suivi de \cite[Theorem 1]{Zhu} en prenant en compte \cite[(1.4)]{doiI}. Les équivalences $(1') \Leftrightarrow (1)$ et $(3) \Leftrightarrow (3')$ se montrent facilement et les équivalences $(2) \Leftrightarrow (3)$ et $(2')\Leftrightarrow (3')$ résultent de \cite[Lemma 22]{Zhu}. \\
Nous allons montrer que $(3') \Rightarrow (1')$. Etant donné les hypothèses du lemme, on peut montrer que $A_C$ est de présentation finie sur $C$ et comme par hypothèse $B$ est plat sur $C$, on peut appliquer \cite[10.11]{wisbauer} pour obtenir les isomorphismes suivants :
$$B\square_{A_C} B^* \simeq Com_{A_C}(B, B)  \ et \ B \square_{A_C} (B \otimes_RA)^*\simeq Com_{A_C} (B \otimes_C A', B)$$ où $\square$ dénote le produit cotensoriel et $(-)^*= Hom_C (-, C)$ est le dual sur $C$.
Puisque $(B\otimes \epsilon)$ est une section $C$-linéaire de $\rho_B : B \rightarrow B \otimes_R A$, le morphisme $(B\otimes_C A_C)^* \rightarrow B^*$ est surjectif et puisque nous supposons $B$ $A_C$-coplat, $B\square_{A_C} (B \otimes_RA)^* \rightarrow B \square_{A_C} B^*$ est aussi surjective. C'est ainsi, que par les isomorphismes précédents, le morphisme $Com_{A_C} (B \otimes_C A_C, B)\rightarrow Com_{A_C} (B , B)$ est surjectif, ce qui permet de conclure. 
\end{proof}
\section{Preuve du théorème \ref{thm3}}
\begin{lemme} 
Supposons $G$ fini, plat sur $S$. Les assertions suivantes sont équivalentes :
\begin{enumerate} 
\item Le champ classifiant $[S/G]$ est modéré. 
\item $G$ est linéairement réductif au sens de \cite{AOV}.
\item Pour tout point  géométrique $\bar{x}: Spec \ k \rightarrow S$ de $S$, la fibre $G_{\bar{x}}:= G \times_S Spec \ k$ est linéairement réductive.  
\end{enumerate} 
\end{lemme}
\begin{proof} L'équivalence $(1) \Leftrightarrow (2)$ est immédiate sachant que l'identité $S\rightarrow S$ est dans ce cas un espace de modules grossier. L'équivalence $(1) \Leftrightarrow (3)$ est alors une réécriture de l'équivalence $(a) \Leftrightarrow (b)$ de \cite[Theorem 3.2]{AOV} pour une action triviale. \end{proof}
Après avoir montré que pour tout point $\zeta : Spec (T) \rightarrow X$, $\underline{Aut}_T(p(\zeta)) \simeq \mathcal{I}_G( \zeta )$, où $p$ est le morphisme canonique surjectif $X \rightarrow [X/G]$ et pour tout $\xi : T\rightarrow [X/G]$,  $$\underline{Aut}_T(\xi) := [X/G] \times_{([X/G]\times_S [X/G] , \Delta , \Delta \circ \xi )} Spec(T)$$ 
Prenons un $T$-point $\zeta : Spec \ T \rightarrow X$ au dessus de $S$, comme par hypothèse le groupe d'inertie $\mathcal{I}_G( \zeta )$ est fini et plat, par le lemme précédent, il nous suffit donc de démontrer que $(\mathcal{I}_G( \zeta ))_{\bar{x}}$ est linéairement rédutif, pour tout point $\bar{x}: Spec \ k \rightarrow Spec \ T$, avec $k$ algébriquement clos. Or, on montre facilement que $(\mathcal{I}_G( \zeta ))_{\bar{x}} = \mathcal{I}_G( \zeta \circ \bar{x} )$ où $\zeta \circ \bar{x} : Spec \ k \rightarrow X$ est un point géométrique de $X$. Ainsi, pour tout point géométrique $\bar{x}$, par l'équivalence  $(a) \Leftrightarrow (b)$ de \cite[Theorem 3.2]{AOV}, tous les groupes d'inertie aux points géométriques sont linéairement réductifs au sens de \cite{AOV}, donc $(\mathcal{I}_G( \zeta ))_{\bar{x}}$ est linéairement réductif au sens de \cite{AOV} ce qui montre par le lemme précédent que $\mathcal{I}_G( \zeta )$ l'est aussi. La réciproque est évidente en considérant toujours la même équivalence de \cite[Theorem 3.2]{AOV} .

\small
\bibliography{biblio1} 
\Large
\end{document}